# NONTANGENTIAL AND PROBABILISTIC BOUNDARY BEHAVIOR OF PLURIHARMONIC FUNCTIONS

By Steve Tanner[1]

*Eastern Oregon University*

Let $u$ be a pluriharmonic function on the unit ball in $\mathbb{C}^n$. I consider the relationship between the set of points $L_u$ on the boundary of the ball at which $u$ converges nontangentially and the set of points $\mathcal{L}_u$ at which $u$ converges along conditioned Brownian paths. For harmonic functions $u$ of two variables, the result $L_u \stackrel{\text{a.e.}}{=} \mathcal{L}_u$ has been known for some time, as has a counterexample to the same equality for three variable harmonic functions. I extend the $L_u \stackrel{\text{a.e.}}{=} \mathcal{L}_u$ result to pluriharmonic functions in arbitrary dimensions.

**1. Introduction.**

1.1. *Definitions and notation.* For $r \in (0, \infty)$, let $r\mathbb{B}^n = \{z \in \mathbb{C}^n : \|z\| < r\}$. When there is no confusion as to dimension, I will shorten this notation simply to $r\mathbb{B}$. Further, when referring to the unit ball, I write $\mathbb{B}$ in place of $1\mathbb{B}$.

I use $(Z_t^{(n)}, \mathbb{P}^w, \mathcal{F}_t)$ to denote a complex Brownian motion started at $w \in \mathbb{B} \subset \mathbb{C}^n$. As above, when no confusion may arise, I will shorten this notation simply to $(Z_t, \mathbb{P}^w, \mathcal{F}_t)$.

If $B \subset \mathbb{B}$, define the usual stopping times

$$\tau_B = \inf\{t > 0 : Z_t \in B^c\}$$

and

$$T_B = \inf\{t > 0 : Z_t \in B\}.$$

The main results concern the relationships between two different types of convergence described below.

Received September 2004; revised November 2005.
[1]Supported in part by a Faculty Scholars grant from Eastern Oregon University.
*AMS 2000 subject classifications.* Primary 60J45, 32A40; secondary 60J65.
*Key words and phrases.* Nontangential convergence, pluriharmonic function, conditional process, potential, hitting probability.

DEFINITION 1.1. Let $\theta \in \partial \mathbb{B}$ and let $A \in \mathcal{F}_t$. Define a probability measure $\mathbb{P}_\theta^0$ by the $h$-path transform

$$P_\theta^0(A) = \frac{\mathbb{E}^0[h_\theta(Z_{t \wedge \tau_\mathbb{B}}); A]}{h_\theta(0)},$$

where

$$h_\theta(w) = (1 - \|w\|^2)/\|\theta - w\|^{2n}$$

is the Poisson kernel for $\mathbb{B}$ with a pole at $\theta$ [1].

DEFINITION 1.2. Let $\theta \in \partial \mathbb{B}$. A function $f : \mathbb{B} \to \mathbb{C}$ is said to converge to $L$ along conditioned Brownian paths at $\theta$ if $\lim_{t \to \tau_\mathbb{B}} f(Z_t) = L$ $\mathbb{P}_\theta^0$-a.s. (Note, if the limit exists for a given $\theta$, the Brossard zero–one law [1] guarantees that $L$ is not just a random variable, but actually constant.)

If $\theta \in \partial \mathbb{B}$, define the Stoltz domain $C_{\theta, r}$ with vertex $\theta$ and width $r$ to be the convex hull of the point $\theta$ and the closure of the set $r\mathbb{B}$.

DEFINITION 1.3. Let $\theta \in \partial \mathbb{B}$. A function $f : \mathbb{B} \to \mathbb{C}$ is said to converge nontangentially to $L$ at $\theta$ if $\lim_{z \to \theta, z \in C_{\theta, r}} f(z) = L$ for all $r \in (0, 1)$.

For a function $f : \mathbb{B} \to \mathbb{C}$, define the set of points on the boundary of $\mathbb{B}$ at which this function converges in each of the senses above.

DEFINITION 1.4. If $f : \mathbb{B} \to \mathbb{C}$, define the following sets:

$$L_f = \{\theta \in \partial \mathbb{B} : f \text{ converges nontangentially at } \theta\}$$

and

$$\mathcal{L}_f = \{\theta \in \partial \mathbb{B} : f \text{ converges along conditioned Brownian paths at } \theta\}.$$

1.2. *History and results.*

DEFINITION 1.5. If $A$ and $B$ are Borel sets, we write $A \stackrel{\text{a.e.}}{\subset} B$ if there is a set $C$ of measure zero such that $A \subset B \cup C$.

The a.e. equivalence of the sets $L_u$ and $\mathcal{L}_u$ for harmonic functions on $\mathbb{B} \subset \mathbb{C}^1$ was first studied by Doob [4] using the method of $h$-path transforms. The containment $L_u \stackrel{\text{a.e.}}{\subset} \mathcal{L}_u$ applies to harmonic functions in any dimension (see [1] or [5]), but Burkholder and Gundy produced a counterexample to show that, even for harmonic functions defined in $\mathbb{R}^3$, the containment $\mathcal{L}_u \subset L_u$ need not hold. In fact, one can construct a harmonic function which does not converge nontangentially at any point, but which converges along conditioned Brownian paths on a set of positive measure [3]. In higher dimensions, relations between $\mathcal{L}_u$ and $L_u$ for a more restrictive set of functions were studied by Durett, Chen and Ma [6], who showed a similar containment may be recovered for functions harmonic with respect to a suitably chosen differential operator, even when nontangential convergence is replaced with convergence in larger admissible regions.

My results show that the $\mathcal{L}_u \subset L_u$ holds if one requires the function $u$ to be not only harmonic, but pluriharmonic—that is, the real or imaginary part of a holomorphic function on $\mathbb{C}^n$. Specifically, I prove the following:

THEOREM 1.1. *If $u : \mathbb{B} \subset \mathbb{C}^n \to \mathbb{R}$ is pluriharmonic, then $\mathcal{L}_u \subset L_u$.*

Since all pluriharmonic functions are harmonic, the following corollary follows from the theorem and the older results cited above.

COROLLARY 1.1. *If $u : \mathbb{B} \subset \mathbb{C}^n \to \mathbb{R}$ is pluriharmonic, then $\mathcal{L}_u \stackrel{a.e.}{=} L_u$.*

## 2. Brownian convergence implies nontangential convergence.

2.1. *Motivation and layout.* The Burkholder–Gundy counterexample in $\mathbb{R}^3$ is produced by constructing what the authors refer to as "a bed of nails." It relies on the fact that a harmonic function may be constructed so that it has a large modulus only on a series "spikes" getting progressively thinner near a boundary. This prevents the function from converging nontangentially, but since the "spikes" may be chosen small enough that Brownian motion will not hit them too often, the conditioned Brownian paths will eventually fail to intersect the offending sets and the function will converge along such paths. Pluriharmonicity is a stronger restriction which requires functions to be harmonic when restricted to any complex line [7] and prevents such a construction.

The idea behind the proof is to show that if a pluriharmonic function is large at a point, then it is large on a set with capacity bounded away from zero. With that in mind, the remainder of the Section 2 is as follows. In Section 2.2 it is proven that, around any point in a Stoltz cone, balls of a certain relative size may be hit by Brownian motion with some minimal probability. A mapping $\Phi$ that is useful for indexing various complex lines through the origin is defined and some of its properties are noted in Section 2.3. A technical lemma concerning the continuity and boundedness of the Green potentials of a certain class of measures is proven in Section 2.4. The $\Phi$ function from above is used in Section 2.5 to create a measure in the class dealt with from the previous section. Finally, in Section 2.6 this measure is used to prove Theorem 1.1.

2.2. *Probability estimates for balls.*

PROPOSITION 2.1. *Let $\theta \in \partial \mathbb{B}$ and let $C_{\theta,r}$ be a Stoltz domain. There is a constant $c_{2.1} > 0$ such that, if $z \in C_{\theta,r}$ and $B_z = B(z, (1 - \|z\|)/8)$, then $\mathbb{P}^0_\theta(T_{B_z} < \tau_\mathbb{B}) \geq c_{2.1}$.*

PROOF. There is a real number $s_r \in (0,1)$ so that if $z \in C_{\theta,r}$ and $w \in B_z$, then $w \in C_{\theta,s_r}$, and therefore, there is a positive constant $c_1$ so that $\|\theta - w\| \leq c_1(1 - \|z\|)$. Let $C = \partial(\|z\|\mathbb{B}) \cap B_z$. Note that there is a positive constant $c_2$ so that $\sigma(C) \geq c_2(1 - \|z\|)^{2n-1}$ if $\sigma$ is a normalized surface measure on $\|z\|\mathbb{B}$. Let $T_z = \inf\{t : \|Z_t\| = \|z\|\}$. By the spherical symmetry of Brownian motion, there is a constant $c_3$ so that $\mathbb{P}^0(Z_{T_z} \in C) \geq c_3 \sigma(C)$:

$$\mathbb{P}^0_\theta(T_{B_z} < \tau_\mathbb{B}) \geq \mathbb{P}^0_\theta(Z_{T_z \in C})$$
$$= \mathbb{E}^0\left[\frac{h_\theta(Z_{T_z})}{h_\theta(0)}; Z_{T_z} \in C\right]$$
$$\geq \mathbb{E}^0\left[\frac{(1 - \|Z_{T_z}\|^2)}{\|\theta - Z_{T_z}\|^{2n}}; Z_{T_z} \in C\right]$$
$$\geq c_1 \mathbb{P}^0(Z_{T_z} \in C)\frac{1 - \|z\|}{(1 - \|z\|)^{2n}}$$
$$\geq c_1 c_3 \sigma(C)\frac{1 - \|z\|}{(1 - \|z\|)^{2n}}$$
$$\geq c_1 c_2 c_3 (1 - \|z\|)^{2n-1}\frac{1 - \|z\|}{(1 - \|z\|)^{2n}}$$
$$= c_{2.1} > 0. \qquad \square$$

2.3. *The $\Phi$ mapping.*

DEFINITION 2.1. Define $\Phi : [0, 2\pi) \times (0,1) \times (\mathbb{C}^{n-1}\setminus\{0\}) \to \mathbb{C}^n$ by

$$\Phi(\theta, r, r(2)e^{i\theta(2)}, \ldots, r(n)e^{i\theta(n)})$$
$$= \left[\frac{r}{\sqrt{1 + \sum_{j=2}^n r(j)^2}}\right](e^{i\theta}, r(2)e^{i(\theta - \theta(2))}, \ldots, r(n)e^{i(\theta - \theta(n))})$$

and note that $\Phi$ is a one-to-one function (we assume $r_j > 0$ and $\theta_j \in [0, 2\pi)$ for $j = 2, \ldots, n$).

Let $A = \{w \in \mathbb{C} : 1/3 \leq \|w\| \leq 2/3\}$ and let $A^k$ denote the cross product of $A$ with itself $k$ times. Let $\Delta = [1/8, 1/4] \times A^{n-1}$. Note that $\Phi$ maps

$[0, 2\pi) \times \Delta$ into $\mathcal{A} = \{w \in \mathbb{C}^n : 1/8 \leq \|w\| \leq 1/4\}$. Also note that if $(w_1, \ldots, w_n) = \Phi(\theta, r, z_2, \ldots, z_n)$, then

$$0 < \frac{1/8}{\sqrt{1 + (n-1)(2/3)^2}} \leq \|w_1\|$$

$$\leq \frac{1/4}{\sqrt{1 + (n-1)(1/3)^2}} < \frac{1}{4},$$

and for $j = 2, \ldots, n$,

$$0 < \frac{(1/3)(1/8)}{\sqrt{1 + (n-1)(2/3)^2}} \leq \|w_j\|$$

$$\leq \frac{(2/3)(1/4)}{\sqrt{1 + (n-1)(1/3)^2}} < \frac{1}{4}.$$

In particular, there is a positive constant $b_-$ depending only on the dimension $n$ such that if $(w_1, \ldots, w_n) \in \Phi([0, 2\pi) \times \Delta)$, then $0 < b_- \leq \|w_j\| \leq 1/4$ for all $j = 1, \ldots, n$. Without loss of generality, take $b_- < 1/(16n^{1/2})$.

Let $\Delta'$ be such that

$$[0, 2\pi) \times \Delta' = \Phi^{-1}\left(\left\{w \in \mathbb{C}^n : \frac{1}{8} - \frac{b_-}{2} \leq \|w\| \leq \frac{1}{4} + \frac{b_-}{2};\right.\right.$$

$$\left.\left.\|w_j\| \geq \frac{b_-}{2} \text{ for } j = 1, \ldots, n\right\}\right)$$

and let $\Delta''$ be such that

$$[0, 2\pi) \times \Delta'' = \Phi^{-1}\left(\left\{w \in \mathbb{C}^n : \frac{1}{8} - \frac{b_-}{4} \leq \|w\| \leq \frac{1}{4} + \frac{b_-}{4};\right.\right.$$

$$\left.\left.\|w_j\| \geq \frac{3b_-}{4} \text{ for } j = 1, \ldots, n\right\}\right).$$

The choice of $b_-$ guarantees both $\Delta'$ and $\Delta''$ are nonempty. Further, since $b_- < 1/16$, there are constants $0 < c_- < c_+ < \infty$ such that

$$\Delta \subset \Delta'' \subset \Delta' \subset [\tfrac{1}{16}, \tfrac{5}{16}] \times (\{w \in \mathbb{C} : c_- \leq \|w\| \leq c_+\})^{n-1} \overset{\mathrm{df}}{=} \Delta_0.$$

2.4. *Potentials.* Throughout Section 2.4 let

$$u(w, z) = \frac{c_{2n}}{\|w - z\|^{2n-2}} \qquad \text{where } c_{2n} = \frac{\Gamma(n-1)}{(2\pi)^n}$$

denote the Newtonian potential in $\mathbb{R}^{2n}$,

$$g(w, z) = u(w, z) - \mathbb{E}^w u(X_{\tau_\mathbb{B}}, z)$$

denote the Green function for $\mathbb{B}$, and

$$H(w) = H_\mu(w) = \int_\mathbb{B} g(w,z)\,d\mu(z)$$

denote the Green potential for measure $\mu$.

LEMMA 2.1. *Let $\mu$ be a measure supported on a set $S \subset \mathbb{B}$. Then $H$ is harmonic on $\mathbb{B}\setminus S$ and $\sup_{w\in\mathbb{B}}\{H(w)\} = \sup_{w\in S}\{H(w)\}$.*

PROOF. See [8], Theorem 1.4 on page 161 and Theorem 1.8 on page 163. □

DEFINITION 2.2. If $s > 0, s > \gamma > 0$, and $\alpha \in [0, 2\pi)$, define the truncated wedge

$$W_{s,\alpha,\gamma} = \{re^{i\theta} \in \mathbb{C} : \|r-s\| \leq \gamma; \|\theta - \alpha\| \leq \gamma\}.$$

The following lemma is modeled on one by Bass and Khoshnevisan (see [2], Proposition 2.7).

LEMMA 2.2. *Let $0 < r_0 < 1/16$. Suppose $\mu$ is a measure on $\mathbb{C}^n$ supported in $\mathcal{A} = \{w \in \mathbb{C}^n : 1/8 \leq \|w\| \leq 1/4\}$ such that $a_- \leq \mu(\mathbb{C}^n) \leq a_+$ for some constants $0 < a_- < a_+$ and satisfying the condition*

$$\sup_{(r(1)e^{i\theta(1)},\ldots,r(n)e^{i\theta(n)})\in\mathcal{A}} \mu(W_{r(1),\gamma,\theta(1),\gamma} \times \cdots \times W_{r(n),\gamma,\theta(n),\gamma}) \leq k\gamma^{2n-1}$$

*for some constant $k > 0$ and all $\gamma < r_0$. Then $H$ is continuous and is bounded on $\mathbb{B}$ by a constant $M_{2.2}$ depending only on $a_-, a_+, k, r_0$ and $n$.*

PROOF. First I show continuity. Let $w \in (r(1)e^{i\theta(1)},\ldots,r(n)e^{i\theta(n)}) \in \mathbb{B}$. If $w \notin \mathcal{A}$, then $H$ is continuous at $w$ since $g(w,z)$ is continuous when $z$ is in the support of $\mu$. Suppose $w \in \mathcal{A}$. Let $\gamma < r_0$, let $W_\gamma = W_{r(1),\gamma,\theta(1),\gamma} \times \cdots \times W_{r(n),\gamma,\theta(n),\gamma}$, and let $z \in W_{\gamma/4}$. The condition stated in the lemma guarantees that $\mu(\{w\}) = 0$. So,

$$\begin{aligned}
|H(w) - H(z)| &= c\left| \int_\mathbb{B} \frac{1}{\|w-\zeta\|^{2n-2}}\,d\mu(\zeta) - \int_\mathbb{B} \frac{1}{\|z-\zeta\|^{2n-2}}\,d\mu(\zeta) \right| \\
&= c\left| \int_{W_{\gamma/2}\cap\mathcal{A}} \frac{1}{\|w-\zeta\|^{2n-2}} - \frac{1}{\|z-\zeta\|^{2n-2}}\,d\mu(\zeta) \right| \\
&\quad + c\left| \int_{\mathcal{A}\setminus W_{\gamma/2}} \frac{1}{\|w-\zeta\|^{2n-2}} - \frac{1}{\|z-\zeta\|^{2n-2}}\,d\mu(\zeta) \right|.
\end{aligned}$$

(1)

Let $(r(1)'e^{i\theta(1)'},\ldots,r(n)'e^{i\theta(n)'})$ denote the coordinates of $z$.

Define $W'_\gamma = W_{r(1)',\gamma,\theta(1)',\gamma} \times \cdots \times W_{r(n)',\gamma,\theta(n)',\gamma}$. Note that $W_{\gamma/2} \subset W'_\gamma$ since $z \in W_{\gamma/4}$. Now let $A_n = W'_{\gamma/2^{n-1}} \setminus W'_{\gamma/2^n}$. To estimate the first term of equation (1), note

$$\int_{W'_\gamma} \frac{1}{\|z-\zeta\|^{2n-2}} \, d\mu(\zeta) \leq \sum_{j=1}^{\infty} \int_{A_j} \frac{1}{\|z-\zeta\|^{2n-2}} \, d\mu(\zeta)$$

$$\leq \sum_{j=1}^{\infty} c\left(\frac{2^j}{\gamma}\right)^{2n-2} \mu(A_j)$$

$$\leq \sum_{j=1}^{\infty} c\left(\frac{2^j}{\gamma}\right)^{2n-2} k(\gamma 2^{-j})^{2n-1}$$

$$\leq \sum_{j=1}^{\infty} c\gamma 2^{-j} = c\gamma.$$

Since $w \in W_{\gamma/4}$ as well, the same estimate applies to $w$. Therefore,

$$\left| \int_{W_{\gamma/2} \cap \mathcal{A}} \frac{1}{\|w-\zeta\|^{2n-2}} - \frac{1}{\|z-\zeta\|^{2n-2}} \, d\mu(\zeta) \right|$$

$$\leq \int_{W_\gamma \cap \mathcal{A}} \frac{1}{\|w-\zeta\|^{2n-2}} + \int_{W'_\gamma \cap \mathcal{A}} \frac{1}{\|z-\zeta\|^{2n-2}} \, d\mu(\zeta) \leq 2c\gamma.$$

To estimate the second term of equation (1), note that since $W_{\gamma/4}$ is bounded away from $\mathcal{A} \setminus W_{\gamma/2}$, then there is a constant $c_n$ depending only on the dimension so that $\|y - \zeta\| \geq c_n\gamma$ for all $y \in W_{\gamma/4}$ and all $\zeta \in \mathcal{A} \setminus W_{\gamma/2}$. In particular, $\|w - \zeta\|$ and $\|z - \zeta\|$ are no less than $c_n\gamma$. Let $c_\zeta = \max\{\|w - \zeta\|, \|z - \zeta\|\}$. Then,

$$\left| \int_{\mathcal{A} \setminus W_{\gamma/2}} \frac{1}{\|w-\zeta\|^{2n-2}} \, d\mu(\zeta) - \int_{\mathcal{A} \setminus W_{\gamma/2}} \frac{1}{\|z-\zeta\|^{2n-2}} \, d\mu(\zeta) \right|$$

$$\leq \int_{\mathcal{A} \setminus W_{\gamma/2}} \left| \frac{(\|w-\zeta\| - \|z-\zeta\|)}{\|w-\zeta\|^{2n-2}\|z-\zeta\|^{2n-2}} \right|$$

(2)
$$\times [\|w-\zeta\|^{2n-3} + \cdots + \|z-\zeta\|^{2n-3}] \, d\mu(\zeta)$$

$$\leq \int_{\mathcal{A} \setminus W_{\gamma/2}} \left| \frac{\|w-z\|}{\|w-\zeta\|^{2n-2}\|z-\zeta\|^{2n-2}} \right| (2n-2)(c_\zeta)^{2n-3} \, d\mu(\zeta).$$

Replacing the larger of $\|w-\zeta\|$ and $\|z-\zeta\|$ by $c_\zeta$ and the smaller by $c_n\gamma$ then gives

$$\left| \int_{\mathcal{A} \setminus W_{\gamma/2}} \frac{1}{\|w-\zeta\|^{2n-2}} \, d\mu(\zeta) - \int_{\mathcal{A} \setminus W_{\gamma/2}} \frac{1}{\|z-\zeta\|^{2n-2}} \, d\mu(\zeta) \right|$$

$$\leq (2n-2) \int_{\mathcal{A}\setminus W_{\gamma/2}} \|w-z\| \frac{(c_\zeta)^{2n-3}}{(c_\zeta)^{2n-2}(c_n\gamma)^{2n-2}} \, d\mu(\zeta)$$

(3)
$$\leq (2n-2) \int_{\mathcal{A}\setminus W_{\gamma/2}} \frac{\|w-z\|}{(c_n\gamma)^{2n-1}} \, d\mu(\zeta)$$

$$\leq \frac{(2n-2)\|w-z\|\mu(\mathbb{C}^n)}{(c_n\gamma)^{2n-1}}.$$

Let $\varepsilon > 0$. Choose $\gamma$ small enough to ensure $2c\gamma < \varepsilon/2$ and then choose $\delta > 0$ such that $(2n-2)\|w-z\|\mu(\mathbb{C}^n)/(c_n\gamma)^{2n-1} < \varepsilon/2$ if $\|w-z\| < \delta$. Then equations (1), (2) and (3) prove the continuity of $H$.

The next step is to show boundedness. Define $A_n = W_{1/2^{n-1}}\setminus W_{1/2^n}$ for $n \geq N$, where $N$ is the integer such that $1/2^N \leq r_0 \leq 1/2^{N-1}$:

$$|H(w)| = H(w) = \int_{\mathbb{B}} g(w,z) \, d\mu(z) \leq \int_{\mathcal{A}} u(w,z) \, d\mu(z)$$

$$\leq \int_{\mathcal{A}\setminus W_{1/2^N}} \frac{c}{\|w-z\|^{2n-2}} \, d\mu(z) + \sum_{j=N}^{\infty} \int_{A_j} \frac{c}{\|w-z\|^{2n-2}} \, d\mu(z)$$

$$\leq \frac{ca_+}{r_0^{2n-2}} + \sum_{j=N}^{\infty} \frac{c}{(2^{-(j+1)})^{2n-2}} \mu(A_j)$$

$$\leq M + c \sum_{j=N}^{\infty} 2^{[(j+1)(2n-2)]} 2^{[(j)(1-2n)]}$$

$$= M + c \sum_{j=N}^{\infty} 2^{-j} = M_{2.2} < \infty. \qquad \square$$

2.5. *Pluriharmonic functions.* Next I construct a measure satisfying the requirements of the previous section and supported on a set on which a pluriharmonic function $u$ has large modulus.

DEFINITION 2.3. For a fixed $\theta \in [0, 2\pi)$, define

$$\Psi_\theta(r, w_2, \ldots, w_n) = \Phi(\theta, r, w_2, \ldots, w_n).$$

Note that, for any $\theta$, $\Psi_\theta$ maps $(r, w_2, \ldots, w_n)$ to a point $(z_1, \ldots, z_n)$ that is both on the complex line $z_1 = w_2 z_2/r_2^2 = \cdots = w_n z_n/r_n^2$, and is also a distance $r$ from the origin.

LEMMA 2.3. *Let $z = (r(1)e^{i\theta(1)}, \ldots, r(n)e^{i\theta(n)}) \in \Psi_\theta(\Delta')$. Denote Lebesgue measure on $\Delta_0$ by $m$. Let $\gamma$ be a real number for which $0 < \gamma < b_-/4$ and let*

$\theta \in [0, 2\pi)$. Then there is a constant $c_{2.3}$ depending only on the dimension $n$ such that

$$m(\Psi_\theta^{-1}(W_{r(1),\gamma,\theta(1),\gamma} \times \cdots \times W_{r(n),\gamma,\theta(n),\gamma})) \leq c_{2.3}\gamma^{2n-1}.$$

PROOF. Note that $\|z\|(1 + \sum_{j=2}^n r_j^2)^{1/2} \geq (1/4)/(1 + (n-1)/9)^{1/2} > 0$, so there is a constant $K$ depending only on the dimension $n$ such that

$$\Psi_\theta^{-1}(W_{r(1),\gamma,\theta(1),\gamma} \times \cdots \times W_{r(n),\gamma,\theta(n),\gamma}) \subset [\|z\| - \gamma, \|z\| + \gamma] \\ \times W_{r(2)/K,\gamma,\theta-\theta(2),\gamma} \times \cdots \\ \times W_{r(n)/K,\gamma,\theta-\theta(n),\gamma}.$$

Thus,

$$m(\Psi_\theta^{-1}(W_{r(1),\gamma,\theta(1),\gamma} \times \cdots \times W_{r(n),\gamma,\theta(n),\gamma})) \\ \leq m([\|z\| - \gamma, \|z\| + \gamma] \times W_{r(2)/K,\gamma,\theta-\theta(2),\gamma} \times \cdots \times W_{r(n)/K,\gamma,\theta-\theta(n),\gamma}) \\ \leq c_{2.3}\gamma^{2n-1}. \qquad \square$$

PROPOSITION 2.2. *Suppose $u$ is pluriharmonic on $\mathbb{B}$, $\|u(0)\| > 0$ and $S = \{w \in \text{int}(\mathcal{A}) : \|u(w)\| > \|u(0)\|/2\}$, where, as before, $\mathcal{A} = \{w \in \mathbb{C}^n : 1/8 \leq \|w\| \leq 1/4\}$. Then there is a constant $k > 0$ (which may be taken independent of $u$) and a measure $\mu$ supported in $S$ so that $\mu$ satisfies the conditions in Lemma 2.2.*

PROOF. $u$ is pluriharmonic, so its restriction to any complex line is a harmonic function and thus obeys the maximum modulus principle [7]. Therefore, given any point $(r, w_2, \ldots, w_n) \in \Delta$, there must exist a point $z \in \Phi(\cdot, r, w_2, \ldots, w_n)$ such that $\|u(z)\| \geq \|u(0)\| > \|u(0)\|/2$. In particular, there exists a $\theta_z \in [0, 2\pi)$ such that $\Psi_{\theta_z}(r, w_2, \ldots, w_n) = (r(1)e^{i\theta(1)}, \ldots, r(n)e^{i\theta(n)}) \in S$. $S$ is open, so there exists an open poly-wedge $W_z = W_{r(1),\gamma,\theta(1),\gamma} \times \cdots \times W_{r(n),\gamma,\theta(n),\gamma}$ such that $W \subset S$. By continuity, $V_z = \Psi_{\theta_z}^{-1}(W_z)$ is open. Repeating this for every point in $\Delta$ and using compactness yields a finite collection of open sets $V_{z_1}, \ldots, V_{z_m}$ covering $\Delta$, and corresponding functions $\Psi_{\theta_{z_1}}, \ldots, \Psi_{\theta_{z_m}}$ mapping these sets into $\mathcal{A}$. In order to simplify later notation, we will artificially name the angles $\theta_{z_0} = 0$ and $\theta_{z_{m+1}} = 2\pi$, and we will assume that $0 = \theta_{z_0} < \theta_{z_1} < \theta_{z_2} < \cdots < \theta_{z_m} < \theta_{z_{m+1}} = 2\pi$.

Let $W_1 = V_{z_1} \cap \Delta$ and $W_j = (V_{z_j} \setminus \bigcup_{k=1}^{j-1} V_{z_k}) \cap \Delta$ if $j > 1$. Then define $\phi(r, w_2, \ldots, w_n) = \Psi_{\theta_{z_j}}(r, w_2, \ldots, w_n)$ if $(r, w_2, \ldots, w_n) \in W_j$. Note that $\phi$ is a well-defined function from $\Delta$ into $S$ which is piecewise continuous.

Next define a measure $\mu$ on $\mathbb{B}$ as follows. If $A \subset \mathbb{B}$,

$$\mu(A) = m(\phi^{-1}(A)),$$

where $m$ is a Lebesgue measure on $\Delta$. Note that $\mu$ is supported in $S$.

All that remains to be shown is that $\mu$ satisfies the conditions of Lemma 2.2. Let $w = (r(1)e^{i\theta(1)}, \ldots, r(n)e^{i\theta(n)}) \in \mathbb{B}$ and let $0 < \gamma < b_-/4$.

Let $W = W_{r(1),\gamma,\theta(1),\gamma} \times \cdots \times W_{r(n),\gamma,\theta(n),\gamma}$. To obtain an estimate on $\mu(W)$, note that, unless $w \in \Psi_\theta(\Delta')$, $W$ lies entirely outside of $\Psi_\theta(\Delta)$ and thus outside of the support of $\mu$. Since, by construction, $\mu$ is supported only on the set where $\theta_1 = \theta_{z_j}$ for some $j = 1, \ldots, m$, the set $W_{r(1),\gamma,\theta(1),\gamma}$ decomposes in the following manner:

$$W_{r(1),\gamma,\theta(1),\gamma} = \{re^{i\theta} \in \mathbb{C} : \|r - r(1)\| \leq \gamma; \|\theta - \theta(1)\| \leq \gamma\}$$

$$= \bigcup_{j=0}^{m} \{re^{i\theta} \in \mathbb{C} : \|r - r(1)\| \leq \gamma;$$

$$\theta \in [\theta_{z_j}, \theta_{z_{j+1}}) \cap [\theta(1) - \gamma, \theta(1) + \gamma]\}.$$

Then by Lemma 2.3,

$$\mu(W_{r_1,\gamma,\theta_1,\gamma} \times W_{r_2,\gamma,\theta_2,\gamma} \times \cdots \times W_{r_n,\gamma,\theta_n,\gamma})$$

$$= m(\phi^{-1}(W_{r_1,\gamma,\theta_1,\gamma} \times W_{r_2,\gamma,\theta_2,\gamma} \times \cdots \times W_{r_n,\gamma,\theta_n,\gamma}))$$

$$= m\left(\phi^{-1}\left(\bigcup_{j=0}^{m} \{re^{i\theta} \in \mathbb{C} : \|r - r_1\| \leq \gamma;\right.\right.$$

$$\theta \in [\theta_{z_j}, \theta_{z_{j+1}}) \cap [\theta_1 - \gamma, \theta_1 + \gamma]\}$$

$$\left.\left.\times W_{r_2,\gamma,\theta_2,\gamma} \times \cdots \times W_{r_n,s,\theta_n,s}\right)\right)$$

$$= m\left(\bigcup_{j=1}^{m} \Psi_{\theta_{z_j}}^{-1}(\{re^{i\theta} \in \mathbb{C} : \|r - r_1\| \leq \gamma;\right.$$

$$\theta \in [\theta_{z_j}, \theta_{z_{j+1}}) \cap [\theta_1 - \gamma, \theta_1 + \gamma]\}$$

$$\left.\times W_{r_2,\gamma,\theta_2,\gamma} \times \cdots \times W_{r_n,\gamma,\theta_n,\gamma})\right)$$

$$\leq \sum_{j=1}^{m} (m(\Psi_{\theta_{z_j}}^{-1}(\{re^{i\theta} \in \mathbb{C} : \|r - r_1\| \leq \gamma;$$

$$\theta \in [\theta_{z_j}, \theta_{z_{j+1}}) \cap [\theta_1 - \gamma, \theta_1 + \gamma]\}$$

$$\times W_{r_2,\gamma,\theta_2,\gamma} \times \cdots \times W_{r_n,\gamma,\theta_n,\gamma})))$$

$$\leq k\gamma^{2n-1}. \qquad \square$$

PROPOSITION 2.3. *Suppose $u$ is pluriharmonic on $\mathbb{B}$, $z \in \mathbb{B}$, and $\|u(z)\| > 0$. Let $r = (1 - \|z\|)/2$, and let $S = \{w \in \overline{B(z, r/2)} \setminus B(z, r/4) : \|u(w)\| > \|u(z)\|/2\}$.*

*Then there is a constant $c_{2.3}$ (which depends on the dimension but may be chosen independent of $u$ and of $z$) such that*

$$\mathbb{P}_\theta^\zeta(T_S < \tau_{B(z,3r/2)}) > c_{2.3}$$

*for every $\zeta \in \overline{B(z, r/4)}$.*

PROOF. First suppose $z$ is the origin (and so $r = 1/2$). By Lemma 2.3, there is a measure $\mu$ on $\mathbb{C}$ supported in $S$ such that $m(\Delta) \leq \mu(\mathbb{C}) \leq m(\Delta'')$ and such that $\mu$ satisfies the conditions of Lemma 2.2. Therefore, by Lemmas 2.1 and 2.2, $H(\zeta)$ is a function bounded by a constant $M$ (independent of $u$) such that $H$ is harmonic on $\mathbb{B} \backslash S$.

Note that if $\|w\| = 3/4$, it follows that $H(w) = \int_\mathbb{B} g(w,z)\,d\mu(z) \leq c_d/(1/2)^{d-2} = m_d$.

Further, if $w \leq 1/8$, it follows that $H(w) \geq c_d/(3/8)^{d-2} = M_d > m_d$.

Let $T = \inf\{t : H(Z_t) = m_d\} \wedge T_S$ and note that, since $H$ is continuous, $T < \tau_{(3/4)\mathbb{B}}$. Let $Z_t$ be a Brownian motion on $\mathbb{B}$. Since $H$ is harmonic on $\mathbb{B} \backslash S$, $W_t = H(Z_{t \wedge T})$ is a continuous martingale.

Let $M = M_{2.2}$, the upper bound on $H$ which was found earlier. Note that $Z_T \in S$ if and only if $W_T > M_d$, so $\mathbb{P}^\zeta(T_S < \tau_{(3/4)\mathbb{B}}) = \mathbb{P}^\zeta(W_T \geq M_d)$.

Either $W_T = m_d$ or $W_T \in S$. So if $\zeta \in (1/8)\mathbb{B}$, then

$$\mathbb{P}^\zeta(W_T = m_d) + \mathbb{P}^\zeta(W_T \geq M_d) = 1$$

and

$$m_d \mathbb{P}^\zeta(W_T = m_d) + M \mathbb{P}^\zeta(W_T \geq M_d) \geq \mathbb{E}^\zeta[W_T] \geq M_d.$$

Combining these two equations yields

$$\mathbb{P}^\zeta(W_T \geq M_d) \geq \frac{M_d - m_d}{M - m_d} = c > 0.$$

Next suppose $z$ is an arbitrary point in $\mathbb{B}$. Let $\tilde{u}(w) = u(z - rw)$. So $\|\tilde{u}(0)\| = \|u(z)\| > 0$ and $\tilde{u}$ is pluriharmonic in $\mathbb{B}$.

Let $S_{\tilde{u}} = \{w \in \overline{(1/4)\mathbb{B}} \backslash (1/8)\mathbb{B} : \|\tilde{u}(w)\| > \|\tilde{u}(0)\|/2\}$. Let $\tilde{Z}_t$ be a Brownian motion and let $\tilde{T}_{S_{\tilde{u}}}$ and $\tilde{\tau}_\mathbb{B}$ be the corresponding stopping times. Then, from above, $P^{\tilde{\zeta}}(\tilde{T}_{S_{\tilde{u}}} < \tilde{\tau}_\mathbb{B}) > c_{2.3}$ for every $\tilde{\zeta} \in (1/8)\mathbb{B}$.

Now define $Z_t = z - r\tilde{Z}_{t/r^2}$, which by scaling and translation invariance is a Brownian motion. If $\zeta \in \overline{B(z, r/4)}$, then $\mathbb{P}^\zeta(T_S < \tau_{B(z,3r/2)}) = \mathbb{P}^{\tilde{\zeta}}(\tilde{T}_{S_{\tilde{u}}} < \tilde{\tau}_\mathbb{B}) > c > 0$.

Next I obtain a similar result for the conditional probabilities. By the Harnack principle, there is a constant $c$ such that, for $w \in \overline{B(z, 3r/2)}$ and, for $\zeta \in \overline{B(z, r/4)}$,

$$\frac{h_\theta(w)}{h_\theta(\zeta)} \geq c > 0.$$

So,

$$\mathbb{P}_\theta^\zeta(T_S < \tau_{B(z,3r/2)}) = E^\zeta\left[\frac{h_\theta(Z_{T_S \wedge \tau_{B(z,3r/2)}})}{h_\theta(\zeta)}; T_S < \tau_{B(z,3r/2)}\right]$$
$$\geq c\mathbb{P}^\zeta(T_S < \tau_{B(z,3r/2)})$$
$$\geq c_{2.3} > 0. \qquad \square$$

2.6. *Convergence result.*

PROOF OF THEOREM 1.1. Suppose the containment is false. Let $\theta \in \mathcal{L}_u \setminus L_u$. Without loss of generality, assume that $u$ converges to 0 along conditioned Brownian paths at $\theta$. Since $\theta$ is not in $L_u$, there is a Stoltz cone $C_{r,\theta}$ and a sequence of points $z_n \in C_{r,\theta}$ converging to $\theta$ such that $\|u(z_n)\| \geq \varepsilon$ for some $\varepsilon > 0$.

Let $A_n = \{w : (1-\|z_n\|)/8 \leq \|w-z_n\| \leq (1-\|z_n\|)/4\}$, let $B_n = \{w : \|w-z_n\| \leq (1-\|z_n\|)/8\}$, and let $S_n = \{w \in B_n : \|u(w)\| > \|u(z_n)\|/2\}$. Then, by the strong Markov property, Propositions 2.1, 2.3 and the fact that Brownian motion must eventually exit $\mathbb{B}$, we have

$$\mathbb{P}_\theta^0(T_{S_n} < \tau_\mathbb{B}) = \mathbb{E}_\theta^0[\mathbb{P}_\theta^0(T_{S_n} < \tau_\mathbb{B})]$$
$$= \mathbb{E}_\theta^0[\mathbb{P}_\theta^{T_{B_n'}}(T_{S_n} < \tau_\mathbb{B}); T_{B_n} < \tau_\mathbb{B}]$$
$$\geq c_{2.3}\mathbb{P}_\theta^0(T_{B_n} < \tau_\mathbb{B})$$
$$\geq c_{2.3}c_{2.1} > 0.$$

Thus, $\mathbb{P}_\theta^0(T_{S_n} < \tau_\mathbb{B} \text{ i.o.}) > 0$, and so by the Brossard zero–one law, $\mathbb{P}_\theta^0(T_{S_n} < \tau_\mathbb{B} \text{ i.o.}) = 1$. This contradicts the assumption that $\lim_{t \to \tau_\mathbb{B}} f(Z_t) = 0$. $\square$

**Acknowledgment.** The results that follow are based on work done for my Ph.D. thesis and I would like to thank my thesis advisor, Richard Bass, for his help and encouragement in that effort.

MATHEMATICS PROGRAM
EASTERN OREGON UNIVERSITY
LAGRANDE, OREGON 97850
USA
E-MAIL: stanner@eou.edu
URL: www2.eou.edu/~stanner